\documentclass[11pt]{amsart}
\usepackage{amsmath,mathrsfs}

\begin{document}
\theoremstyle{theorem}
\newtheorem{thm}{Theorem}
\newtheorem{lem}{Lemma}
\newtheorem{prop}{Proposition}
\newtheorem{cor}{Corollary}
\theoremstyle{remark}
\newtheorem{rem}{Remark}
\newtheorem*{ack}{Acknowledgement}

\def\Z{{\mathbb Z}} 
\def\C{{\mathbb C}} 
\def\R{{\mathbb R}} 

\def\A{{\mathfrak A}}
\def\B{{\mathcal B}}
\def\H{{\mathcal H}}

\def\L{{\mathfrak L}}
\def\LH{{\mathfrak L(\H)}}
\def\ran{\mathrm{ran}}

\mathchardef\ordinarycolon\mathcode`\:
\def\vcentcolon{\mathrel{\mathop\ordinarycolon}}
\providecommand*\coloneqq{\mathrel{\vcentcolon\mkern-1.2mu}=}

\title{About the proof of the Fredholm Alternative theorems}
\author{Ali Reza Khatoon Abadi }
\address{Ali Reza Khatoon Abadi:Department of  Mathematics,Faculty of science,Islamic azad
university,TEHRAN west Branch}
\author{H.R.Rezazadeh}
\address{H.R.Rezazadeh:Department of  Mathematics,Faculty of
science,Islamic azad university,Karaj branch,Karaj,iran}
\date{\today}
\begin{abstract}
In this short paper we review and extract some features of the
Fredholm Alternative problem .
\end{abstract}

\maketitle
\section{Preliminaries}
In mathematics, the Fredholm alternative, named after Ivar Fredholm,
is one of Fredholm's theorems and is a result in Fredholm theory. It
may be expressed in several ways, as a theorem of linear algebra, a
theorem of integral equations, or as a theorem on Fredholm
operators. Part of the result states that a non-zero complex number
in the spectrum of a compact operator is an
eigenvalue.\cite{Fredholm}
 \subsection{ Linear algebra}
If V is an n-dimensional vector space and $T$ is a linear
transformation, then exactly one of the following holds \\
1-For each vector v in V there is a vector u in V so that $T(u) =
v$. In other words: T is surjective (and so also bijective, since V
is finite-dimensional).\\
2-$dim(Ker(T))>0$\\
A more elementary formulation, in terms of matrices, is as follows.
Given an $m×n$ matrix A and a $m×1$ column vector b, exactly one of
the following must hold
\\
Either: $A x = b$ has a solution x Or: $A^{T }y = 0$ has a solution
y with $y^{T}b \neq 0$.\\
In other words, $A x = b$ has a solution  if and only if for any
$A^{T }y = 0$,$y^{T}b = 0$.
\subsection{ Integral equations}
Let K(x,y) be an integral kernel, and consider the homogeneous
equation, the Fredholm integral equation

\begin{eqnarray}
\lambda\phi(x)-\int_{a}^{b}K(x,y)\phi(y)dy=f(x)
\end{eqnarray}
 ,The Fredholm alternative states that, for any non-zero fixed
complex number , either the Homogenous equation has a non-trivial
solution, or the inhomogeneous equation has a solution for all f(x).
A sufficient condition for this theorem to hold is for K(x,y) to be
square integrable on the rectangle $[a,b]\times[a,b]$ (where a
and/or b may be minus or plus infinity).

\subsection{Functional analysis}
Results on the Fredholm operator generalize these results to vector
spaces of infinite dimensions, Banach spaces.
\subsection{Correspondence}
Loosely speaking, the correspondence between the linear algebra
version, and the integral equation version, is as follows: Let $T
=\lambda  - K$ or, in index notation,

$T(x,y) =\lambda \delta(x - y) - K(x,y) $ with $\delta(x - y)$ the
Dirac delta function. Here, T can be seen to be an linear operator
acting on a Banach space V of functions f(x), so that $T:V\mapsto V$
is given by $\phi\mapsto \phi$ with $\psi$ given by
$\psi(x)=\int_{a}^{b}K(x,y)\phi(y)dy$

In this language, the integral equation alternatives are seen to
correspond to the linear algebra alternatives
\subsection{Alternative}
In more precise terms, the Fredholm alternative only applies when K
is a compact operator. From Fredholm theory, smooth integral kernels
are compact operators. The Fredholm alternative may be restated in
the following form: a nonzero ? is either an eigenvalue of K, or it
lies in the domain of the resolvent\cite{Ramm},\cite{Khvedelidze}.

We recall that a bounded linear operator $T:V \to V$ on a Banach
space $V$ is called a {\em compact} operator if $T$ maps the closed
unit ball of $V$ to a relatively compact subset of $V$. The
following is a basic result about compact operators known as the
{\em Fredholm Alternative}.  Let $T$ be a compact operator on a
Banach space
 $V$ and let $\mathfrak{I}$ denote the identity operator on $V$. Then either
the operator $\mathfrak{I} - T$ is invertible (i.e.\ has a bounded
inverse), or
 there exists a nonzero vector $\zeta \in V$, $\zeta \ne 0$, such that $T\zeta = \zeta$.
\section{ Saddle Point's version of Fredholm alternative theorem} Theorem : Let the operator $A: H
\mapsto H$ be linear, compact and self-adjoint on the separable
Hilbert space H. Let $f\in H$, $\lambda\neq0 \in R$ . Then the
problem
\begin{eqnarray}
\lambda u-A u=f
\end{eqnarray}
 has a solution iff
\begin{eqnarray}
(u,f)=0
\end{eqnarray}
for all $u\in H$. There is a beautifull theorem about it\\
 Theorem
: Let the operator $A: H \mapsto H$ be a linear, compact, symmetric
and positive on the separable Hilbert space H and let $\lambda\neq0
\in R$, $\lambda_{1}>\lambda>0$. Then the Fredholm alternative for
the operator A and the equation (2) is a consequence of the Saddle
Point Theorem. For a good proof see the refrence\cite{TOMICZEK}.In
this reference it the author has been proved that the Fredholm
alternative for such an operator is a consequence of the Saddle
Point Theorem.Also see the reference [4] in it\cite{4}.
\section{Unbounded Fredholm Operators}
One can study he topology of the space of all (generally unbounded)
self-adjoint Fredholm operators, and to put the notion of spectral
flow for continuous paths of such operators on a firm mathematical
footing .Also there is a parallel method due to Nicolaescu's
approach which requires the continuity of the Riesz map and to
achieve that additional properties of the families of boundary
problems\cite{Nicolaescu}.Operator curves on manifolds with boundary
In low-dimensional topology and quantum field theory, various
examples of operator curves appear which take their departure in a
symmetric elliptic differential operator of first order (usually an
operator of Dirac type) on a fixed compact Riemannian smooth
manifold M with boundary $\sum$. Posing a suitable well-posed
boundary value problem provides for a nicely spaced discrete
spectrum near 0. Then, varying the coefficients of the differential
operator and the imposed boundary condition suggests the use of the
powerful topological concept of spectral flow. We can show under
which conditions the curves of the induced self-adjoint
$L^{2}$-extensions become continuous curves in  the gap topology
such that their spectral flow is well defined and truly homotopy
invariant.But we remind it to another work.
\section{Summary} We review some aspects of the Fredholm
theorem .We show that there is a relation between topology and this
theorem.

\begin{ack} The authors gratefully acknowledge support from the Islamic Azad University of KARAJ and TEHRAN (West).
\end{ack}


\begin{thebibliography}{35}
\bibitem{Fredholm}
E.I. Fredholm, "Sur une classe d'equations fonctionnelles", Acta
Math. , 27 (1903) pp. 365–390.
\bibitem{Ramm}
A. G. Ramm, "A Simple Proof of the Fredholm Alternative and a
Characterization of the Fredholm Operators", American Mathematical
Monthly, 108 (2001) p. 855
\bibitem{Khvedelidze}
Khvedelidze, B.V. (2001), "Fredholm theorem for integral equations",
in Hazewinkel, Michiel, Encyclopaedia of Mathematics, Springer, ISBN
978-1556080104
\bibitem{TOMICZEK}
PETR TOMICZEK, Proceedings of Equadiff-11 2005, pp. 491–497 ISBN
978-80-227-2624-5
\bibitem{4} P. Tomiczek, The generalization of the
Landesman-Lazer conditon, EJDE, 2001(04) (2001), 1–11.
\bibitem{Nicolaescu}
L. Nicolaescu: On the space of Fredholm operators.
math.DG/0005089 (2000)
\bibitem[1.]{1}
Kantorovich, L., Akilov, G.,  Functional analysis in normed spaces,
Macmillan, New York, 1964

\bibitem{2}
Rudin, W.,  Functional analysis, McGraw Hill, New York, 1973



\end{thebibliography}
\end{document}